\pdfoutput=1
\documentclass[journal]{IEEEtran}
\usepackage{cite}
\usepackage{amsmath,amssymb,amsfonts}
\usepackage{algorithmic}
\usepackage{array}
\usepackage[caption=false,font=normalsize,labelfont=sf,textfont=sf]{subfig}
\usepackage{graphicx}
\usepackage{textcomp}
\usepackage{xurl}
\usepackage{booktabs}
\usepackage{xcolor}
\usepackage{hyperref}
\usepackage{multirow}
\usepackage{makecell}
\usepackage{siunitx}
\usepackage{balance}
\usepackage{float}
\def\BibTeX{{\rm B\kern-.05em{\sc i\kern-.025em b}\kern-.08em
    T\kern-.1667em\lower.7ex\hbox{E}\kern-.125emX}}

\usepackage{scalerel}
\usepackage{orcidlink}

\allowdisplaybreaks

\begin{document}

\title{An Optimized Expansion Auction Integrating Energy, Flexibility, and Network Constraints}

\author{Rafael Benchimol Klausner$^{\orcidlink{0009-0007-9305-5586}}$,
Alexandre Street$^{\orcidlink{0000-0003-1569-030X}}$,~\IEEEmembership{Senior Member,~IEEE}, \\
Joaquim Dias Garcia$^{\orcidlink{0000-0002-7721-8564}}$,~\IEEEmembership{Senior Member,~IEEE},
and~Luiz Barroso$^{\orcidlink{0000-0003-3336-6896}}$,~\IEEEmembership{Fellow,~IEEE}
\thanks{R. B. Klausner is with PSR and LAMPS at PUC-Rio, Rio de Janeiro, Brazil (e-mail: rafabench@psr-inc.com).}
\thanks{A. Street is with the Department of Electrical Engineering at the Pontifical Catholic University of Rio de Janeiro, Brazil (e-mail: street@puc-rio.br), and is a Tinker Visiting Professor at Stanford University’s Center for Latin American Studies (CLAS) and Program on Energy and Sustainable Development (PESD).}
\thanks{J. D. Garcia is with PSR, Rio de Janeiro, Brazil (e-mail: joaquim@psr-inc.com).}
\thanks{L. Barroso is with PSR, Rio de Janeiro, Brazil (e-mail: luiz@psr-inc.com).}
\thanks{R. B. Klausner was partially supported by the Coordenação de Aperfeiçoamento de Pessoal de Nível Superior - Brasil (CAPES) - Finance Code 001. A. Street was also partially supported by FAPERJ and CNPq.}
}

\markboth{Klausner \MakeLowercase{\textit{et al.}}: An Optimized Expansion Auction}%
{Klausner \MakeLowercase{\textit{et al.}}: An Optimized Expansion Auction}

\maketitle

\begin{abstract}
Auction mechanisms are widely used to support investments in generation investments. In the US, capacity auctions serve as supply-adequacy instruments to address the missing-money problem, whereas in Latin America, energy and capacity auctions procure long-term contracts to support generation expansion. In both settings, accreditation metrics, such as Firm Energy (FE) and Effective Load Carrying Capability (ELCC), are used to approximate the reliability contribution of selected units and winning projects are selected based on a financial offer. This approach fails to recognize the actual contribution of winners to the multiple services they provide to the system. This limitation becomes more significant with high shares of renewables, as the value of each resource increasingly depends on network and technological constraints, as well as on complementarities within the generation portfolio. Inspired by short-term auctions, such as day-ahead markets, this paper proposes a joint optimization-based auction for generation and transmission investments that internalizes these interactions, capturing the intrinsic value of generation candidates without requiring the computation of FE or ELCC metrics. The auctioneer determines the winning portfolio through an expansion-planning model that explicitly represents network constraints, operational interactions, reliability standards, and multiple scenarios. A case study of the Brazilian system indicates lower system costs, lower risk, and a more diverse portfolio relative to the stylized decoupled procurement approaches currently deployed in many systems.
\end{abstract}


\section{Introduction}

\IEEEPARstart{E}{lectricity} systems worldwide are undergoing structural changes that are increasing pressure on the pace of capacity expansion and on the criteria used to evaluate individual generator contributions to system adequacy. First, the large-scale integration of wind and solar generation has introduced pronounced intraday variability, challenging planning frameworks that have historically focused on long-term energy adequacy~\cite{BritoPereira2022adjusting}. This shift has increased the system value of flexibility, such as fast-ramping generation, storage, and demand response, and highlighted the need for expansion criteria that explicitly account for these attributes. Second, load growth has re-emerged as a binding constraint. While sustained demand growth has long characterized developing economies, it is now also evident in developed systems, driven by electrification and the rapid expansion of data centers~\cite{IEA_Electricity2025}. This has revived an interest in support mechanisms to secure generation investments, including long-term energy and capacity auctions.

Long-term contract auctions have become the primary instrument for many countries to ensure that capacity expansion keeps pace with demand growth~\cite{MaurerBarroso2011}, especially in Latin America. Capacity auctions also remain widely adopted worldwide (see \cite{oren2003generation} for a didactic and insightful discussion). Generally, both approaches require assigning a capacity (or energy) credit for each generation unit. If the capacity credit is set too high, the system is less reliable than the auction implies. If it is set too low, the system over-procures capacity to meet reliability standards. This issue mirrors long-standing concerns in the energy economics literature about how standard cost metrics obscure substantive differences across technologies. Joskow~\cite{Joskow2011} shows that conventional levelized cost comparisons implicitly treat all energy as homogeneous, misrepresenting the economic value of intermittent versus dispatchable technologies by ignoring the timing and market value of output. Subsequent work formalizes this insight at the system level: Ueckerdt et al.~\cite{Ueckerdt2013} introduce System LCOE, decomposing integration costs into profile, balancing, and grid components, while Idel~\cite{Idel2022} demonstrates that full-system cost metrics are inherently portfolio-dependent. These contributions highlight a central implication for auction design: adequacy and value cannot be assessed at the level of isolated projects. Integrating such temporal, locational, and operational heterogeneity into capacity valuation is therefore essential for eliciting truthful bids, a central concern of optimal auction design~\cite{Myerson1981}, and delivering reliable capacity at least cost.

A widely adopted capacity credit methodology is the Effective Load Carrying Capability~\cite{NREL2024ELCC}. It measures how much additional demand the system can reliably serve after adding a given resource to the existing fleet. However, achieving adequacy through these mechanisms is not the same as achieving it efficiently. Marshman et al.~\cite{MARSHMAN2020104803} show that revenue insufficiency in high-renewable systems results primarily from non-optimal fleet composition, specifically not allowing the retirement of thermal generators, not variability-induced market failure. Wolak~\cite{Wolak2021} argues more broadly that long-term resource adequacy under high renewables requires procurement mechanisms that explicitly account for intermittency and portfolio interactions. Gonz\'alez-Romero et al.~\cite{GonzalezRomero2020} show that coordinated generation-transmission investment captures technological synergies that decoupled procurement overlooks.

Latin America is a good example of a region that has been relying on this recipe by carrying out long-term auctions for more than 20 years. Brazil is the leading country in this task and in order to foster new generation investments the country carries out exclusive auctions for new capapcity, known as New Energy Auctions (\textbf{NEA}s). In this framework, the NEAs rely on (1) firm energy certificates, which aim to quantify the firm energy contribution of each resource~\cite{Faria2009firmEnergy}, and (2) a cost-and-benefit index, a scoring rule to estimate variable costs of each technology and compare them on a \$ / MWh basis ~\cite{bezerra2010bidding}. Interestingly, both indices are administratively calculated and rely on simulations of the system performed prior to the auction and ignore the impact of transmission and short-term operational constraints. As a result, the ranking of projects in the auction may not reflect their actual economic value in real system operation. 

In March 2026, Brazil cleared its second Capacity Auction (\textbf{CA}), procuring fixed-payment capacity products segmented by region and resource type based on needs identified in prior expansion planning studies using administrative data~\cite{LRCAP_2026_results}. The auction was not technology-neutral: its product definition favored dispatchable thermal capacity, and it ultimately contracted 16.47~GW of thermal generation against only 2.5~GW of hydro. The resulting obligation exceeds R\$~38.9~billion per year (approximately \$7.5~billion at the prevailing exchange rate), equivalent to about \$395{,}000 per MW-year---roughly 3.3 times the clearing price of PJM's most recent annual capacity auction (the 2028/29 delivery year)~\cite{PJM_BRA_2028}. A single capacity mechanism, layered onto the existing market rather than addressing the structural drivers of the emerging capacity needs, thus committed consumers to a multi-decade payment well above comparable international benchmarks.

Within this design, an explicit attempt was made to reward operational flexibility: a deterministic flexibility index, computed from each unit's declared technical parameters, is added to the offered fixed revenue so that faster-ramping units obtain a more favorable evaluation. The adder converts a unit's declared variable cost into an equivalent expected operational cost by assuming a fixed number of annual dispatch hours, identical for every unit, and then uses the declared marginal cost to differentiate this operational expense across the competing thermal plants. While this rudimentary dispatch-cost model does introduce some operational granularity, the assumed number of dispatch hours is fixed administratively rather than derived from system operation, so the resulting ranking reflects a postulated dispatch pattern rather than an optimized one. Despite their differences, these mechanisms share a fundamental limitation: adequacy and flexibility contributions are defined ex ante as administrative parameters, rather than emerging endogenously from system-wide optimization.

In our experience, most current long-term auction designs, especially capacity auctions including those carried out in the US and in Latin America, fail to capture the portfolio-dependent interactions among resources and their overall integrated benefit to the system when selecting generation winners in those auctions. The adequacy value of a given technology depends not only on its intrinsic characteristics, but also on its location, the existing network, and the composition of the contracted generation mix. For example, a wind plant located in the Northeast of part of the system contributes differently to system adequacy than an identical plant in the South, and the value of storage depends critically on the availability of other flexible resources already in the portfolio. These complementarities and diminishing marginal contributions cannot be represented by adequacy metrics computed independently for each project. In the same vein, all estimates of future variable costs of generation projects used to compose scoring rules, such as those of the Brazilian auctions, ignore these complex operational interactions as well.

Beyond resource complementarities, power system expansion also involves coordination challenges among independent investors, as the value of individual projects depends on investments undertaken by other agents~\cite{Corneli2020prism}. Optimization-based procurement can help coordinate these interdependent decisions and reduce information costs while preserving competition among developers.

\subsection{Proposal, Objective, and Contributions}

Planning models are well suited to identifying efficient investment portfolios while internalizing network and technological constraints \cite{Soares2022}, whereas auctions use competition to induce participants to offer their true capacity costs~\cite{Street2011}. Hence, to address the aforementioned issues, this paper proposes a new regulatory mechanism, namely, the Optimized Expansion Auction (\textbf{OEA}), which integrates generation capacity auctions into generation and transmission expansion planning models. The mechanism proposed in this paper does not seek to modify existing contractual and remuneration structures of existing long-term capacity auctions, which warrant further attention and remain open for discussion in a future work. Instead, it focuses on improving the winner-selection process by replacing the current static, \emph{ex ante} predetermined scoring rule with an integrated, optimization-based evaluation of the operational costs and benefits of candidate units under different operating scenarios, allowing the auctioneer to jointly select the winning units as an optimal portfolio.

In our proposed approach, market participants submit offers to the auctioneer composed of their fixed capacity payment, variable costs, and technical parameters, consistent with current qualification practices. The optimization framework simultaneously evaluates each project’s cost and contribution to meeting demand under a given reliability standard across multiple scenarios, explicitly accounting for network constraints, centralized dispatch coordination, and the interaction with the existing generation fleet. Transparency is a core design principle: all input data and scenario assumptions must be auditable, with governance processes that are clearly defined and accessible to all stakeholders. Within this framework, all relevant operational attributes, including energy, ramping capability (often referred to as flexibility), and reserve capabilities, are naturally represented in the auctioneer's (expansion) model. In this context, once offers are submitted, the auction relies on a conventional joint generation and transmission expansion planning model that the planner has developed, tested, and thoroughly validated internally before the auction. The model is then solved to determine the auction winners. This process can also accommodate multiple rounds in which participants can revise their offers, as in an iterative auction~\cite{Street2011}.

The effectiveness of optimization-based procurement mechanisms depends directly on the fidelity of the underlying expansion model. Current implementations, including the one developed in this work, abstract from constraints such as unit commitment, detailed nodal power flow, and higher temporal resolution to maintain tractability in stochastic, multi-temporal frameworks. However, recent advances in tight mixed-integer formulations for unit commitment~\cite{Knueven2020}, convex relaxations of AC power flow~\cite{Low2014a} and scalable decomposition methods~\cite{Soares2022}, continue to push the boundaries of what can be represented within a tractable expansion problem. These modeling advances are amplified by sustained progress in mathematical programming solvers: over the past two decades, combined algorithmic and hardware improvements have yielded LP and MILP speed-ups of roughly two to three orders of magnitude, with algorithmic gains alone making many previously unsolvable instances tractable~\cite{koch2020progress}. 

The central design principle of embedding optimization directly into the market-clearing process is not novel in electricity markets. It already underpins day-ahead market clearing in PJM, CAISO, and ERCOT~\cite{Kirschen2024-pa, Cramton2017}, as well as European power exchanges through the Euphemia algorithm~\cite{nemoeuphemia}. While optimization-based long-term procurement has been previously proposed at a conceptual level~\cite{Corneli2020prism, Joskow2022hybrid}, this work develops and evaluates an auction-clearing mechanism that integrates stochastic generation-transmission co-optimization, quantitatively comparing its outcomes with those of conventional procurement approaches. The proposed auction extends the logic of optimization-based market clearing from short-term operations to long-term investment decisions. Unlike day-ahead dispatch outcomes, which affect system operation for a single day, expansion decisions shape the affordability, reliability, and flexibility of the power system over multiple decades. The proposed approach reverses the traditional sequence of planning followed by capacity auctions. 

More specifically, instead of administratively determining procurement requirements from planning studies based on estimated costs, the optimization model endogenously determines both the quantity and the location of new capacity required to meet the system’s reliability and operational needs. Therefore, the proposed framework advances current industry practice by introducing a new regulatory mechanism that improves existing contracting auction designs for new generation capacity. Specifically, it selects winning projects through a joint generation and transmission expansion planning model that optimizes investment decisions using state-of-the-art modeling approaches (see, e.g.,~\cite{Soares2022} and the references therein). 

Thus, this work aims to develop an alternative auction-clearing mechanism, leading to a new auction design for generation expansion that also accounts for transmission needs. We propose and evaluate the OEA and investigate whether it delivers systematically different, and economically more consistent, outcomes than the conventional decoupled procurement of energy and capacity based on ex-ante resource adequacy metrics. In summary, the paper makes the following contributions:

\begin{itemize}
\item We propose a novel auction design for generation capacity expansion, termed the OEA, which relies on an integrated system-aware procurement process. This framework assesses the systemic value of competing units by simultaneously co-optimizing their operational constraints, network locations, fixed and variable costs, and ramp capacities to form an optimal winning portfolio based on economic and operational offers by bidders. Unlike current decoupled frameworks that rely on administratively defined demand targets segregated by resource type, the OEA internalizes the interactions and synergies between existing resources and new generation and transmission investments within a unified optimization framework. By bridging long-term auction mechanisms with advanced mathematical planning models, this framework leverages the capabilities of modern optimization to ensure strict coherence between the selected portfolio and future market conditions.

\item We quantify the potential gains of the proposed approach relative to conventional decoupled procurement frameworks using detailed and realistic data from the Brazilian power system for the 2034 expansion horizon. Under a common set of investment offers, with identical, technology-typical cost inputs across all mechanisms, we compare how the OEA and conventional decoupled procurement perform in terms of the resulting generation portfolio, total system costs, tail risks, and CO$_2$ emissions.
\end{itemize}

The remainder of the paper is organized as follows. Section~II presents the proposed framework. Section~III describes the case study, while Section~IV reports and discusses the numerical results. Section~V concludes with policy implications and directions for future research.

\section{From Planning Models to Optimized Procurement}

The design of capacity procurement mechanisms involves three interrelated questions: how to define system reliability requirements, how to evaluate individual project contributions, and how to select the winning portfolio. This section addresses each in turn. We first formalize the reliability problem that underlies all procurement frameworks (Section~II-A), then present the proposed mechanism that embeds this formulation directly into the auction clearing process (Section~II-B), and finally contrast it with the current auction models it aims to improve (Section~II-C).

\subsection{Power System Reliability}

A hydrothermal power system generally employs stochastic optimization for both operation and expansion planning. One of the most widely adopted approaches for operational models uses Stochastic Dual Dynamic Programming~\cite{Pereira1991} to calculate the water values under different hydrological conditions and feed short-term operational models. Expansion models generally use L-shaped decomposition~\cite{Gorenstin1993, 9136733, Soares2022} to jointly optimize investment and operational decisions. These models represent hydro-thermal coordination, renewable variability, network constraints, and multi-stage uncertainty.

Resource adequacy, the ability to meet demand reliably over a long-term horizon, can be formulated as a two-stage stochastic program:
\begin{align} \label{Inv_L_shaped}
    \min_{x} \quad & I(x) + \mathbb{E}_{\xi} \left[ O(x,\xi) \right] \\
    \text{s.t.} \quad & \mathcal{R}(x) \leq R^* \label{rel_constraint} \\
    & x \in \mathcal{I}
\end{align}
where $I(x)$ is investment cost, $O(x,\xi)$ is operational cost under scenario $\xi$, $\mathcal{I}$ is the feasible investment set, and $\mathcal{R}(x)$ is a reliability index bounded by threshold $R^*$. Common indices include Loss of Load Expectation and Expected Energy Not Supplied. For energy-constrained systems like Brazil's hydro-dominated grid, energy-based metrics are particularly relevant. In the implementation (Appendix), reliability is enforced implicitly through a penalty on load shedding rather than an explicit constraint on~$\mathcal{R}(x)$.

\subsection{The Optimized Expansion Auction}

The Optimized Expansion Auction rests on two premises: \emph{auctions are efficient mechanisms for price discovery, and expansion planning models are effective at computing optimal portfolios given those prices.} It combines both: agents submit a fixed-revenue bids and technical parameters for their projects, as in conventional auctions, while winner selection uses an integrated expansion and dispatch optimization model rather than simplified metrics. The objective is to select, from submitted bids, the portfolio that minimizes total system cost while satisfying reliability criteria. The auction preserves the existing contract structure, in which developers compete for the market: winning projects receive long-term capacity contracts, as in current auctions. The innovation lies in selection, not in contracting.

In this framework, bidders submit a required fixed annual revenue intended to recover, as much as possible, their annualized investment and fixed operation and maintenance costs, together with a declared variable cost (CVU) reflecting operational expenses. As in any competitive auction, this fixed-revenue bid is a strategic offer and need not coincide with a unit's underlying costs. The winning bidder receives the fixed availability payment throughout the contract horizon, complemented by variable revenue earned in the short-term market based on actual dispatch. Investment recovery therefore depends on the combination of fixed and variable revenue streams rather than on spot revenues alone, and the consistency between these two streams is what determines whether selected projects can recover their costs in operation.

The proposed auction replaces explicit attribute-specific metrics such as Firm Energy, Physical Guarantee or Effective Load Carrying Capability with implicit valuations that emerge from the optimization. Rather than collapsing technical attributes into ex-ante adjustments like the Brazilian Factor~A~\cite{LRCAP_EPE_2023} or the Cost-Benefit Index~\cite{bezerra2010bidding}, it evaluates each project's contribution through dispatch simulation across multiple scenarios, capturing the context-dependent value of flexibility that static adjustments cannot represent. The mechanism is technology-neutral, but not attribute-neutral: hydro, gas, batteries, and demand response compete on equal footing, with contributions evaluated endogenously based on operational performance. Unlike combinatorial auctions, which define multiple products (energy, capacity, flexibility) and combine multi-attribute bids, the OEA internalizes these system requirements directly in the optimization model rather than procuring them as separate products. This shifts the evaluation from individual project cost to total system cost: a project is selected not because it is the cheapest in isolation, but because it produces the greatest reduction in what consumers ultimately pay.

The formulation in Appendix~\ref{mathematical_formulation} details how jointly determining investment and dispatch under uncertainty internalizes externalities that ex-ante capacity credits ignore. The hierarchical time representation (seasons, typical days, and hourly blocks) values technology profiles, hydrological conditions, and intra-day flexibility consistently. Reserve requirements are endogenously met by hydro and battery resources, and transmission flow constraints ensure locational scarcity is reflected in investment selection. Candidates are evaluated jointly with the existing fleet (hydro reservoirs, thermal plants, and installed renewables), so new investments are selected only when they add value beyond what the current system provides, preventing a pitfall of decoupled procurement: contracting new flexibility before unlocking existing flexibility~\cite{IEA_DemandFlex_2025}.

Expansion models exhibit a fundamental primal-dual property: under standard optimality conditions, the revenues derived from shadow prices on energy balance, reserve, and transmission constraints suffice to cover investment costs~\cite{Schweppe1988}. The proposed auction preserves this relationship by construction: project selection and price formation emerge from the same optimization, ensuring winning bids can recover costs through market revenues under the scenarios used for selection. This property holds exactly at the LP relaxation. In the presence of integer investment decisions, shadow prices may not fully support cost recovery for all selected projects, and convex hull pricing or make-whole payments may be required to ensure revenue adequacy~\cite{Schiro2016}. Sequential approaches break this link entirely, Physical Guarantee and capacity payments are computed independently of dispatch, providing no guarantee that projects meeting procurement criteria will earn adequate revenues in operation.

A natural concern is incentive compatibility: because selection depends on portfolio interactions, bidders may shade costs to exploit complementarities. This challenge is inherent to combinatorial environments: Fabra et al.~\cite{Fabra2006} show that strategic behavior depends critically on auction design even in simpler formats. Several factors mitigate  this concern: bid parameters that are subject to engineering-based auditing, a large candidate pool that limits individual market power, and iterative formats such as clock-proxy designs~\cite{Cramton2017} that enables price discovery. A formal game-theoretic analysis, drawing on mechanism design theory~\cite{Myerson1981}, remains an important direction for future work.

The auction can be implemented through single-round sealed bids or multi-round processes. The computational requirements (several hours per run) are compatible with generation project timelines (3--5 years), making auction schedules of 1 to several weeks appropriate.

\subsection{Current Auction Models in Brazil} \label{GF_Brazil}

Brazilian capacity procurement mechanisms can be formalized within a common framework:
\begin{align} \label{Auc_GF}
\min_{x} \quad & \sum_{j \in \mathcal{G}} p_j \phi_j x_j \\
\text{s.t.} \quad & \sum_{j \in \mathcal{G}} x_j \phi_j \geq \hat{d}, \quad x_j \in \{0,1\} \nonumber
\end{align}
where $x_j$ is the selection decision for candidate $j \in \mathcal{G}$, $p_j$ is the unit price bid, and $\phi_j$ quantifies the project's contribution toward the adequacy requirement $\hat{d}$. The two mechanisms differ in how $\phi_j$ and $\hat{d}$ are specified.

In the New Energy Auction (NEA), $\hat{d}$ is the system energy gap and $\phi_j$ is each project's Physical Guarantee— an administratively determined firm energy allocation derived from system-wide reliability studies~\cite{Faria2009firmEnergy}. Bids are ranked by Cost-Benefit Index (annualized cost per unit of Physical Guarantee)~\cite{bezerra2010bidding}. In the Capacity Auction (CA), $\hat{d}$ is  the residual peak capacity gap and $\phi_j = \bar{P}_j (1 - \text{FOR}_j)$. Bids are adjusted with operational attributes such as ramp rates through a multiplicative Factor A and locational adjustments via $\beta$. In both cases, $\phi_j$ is computed independently of the dispatch model used to evaluate the resulting portfolio, and adequacy contributions are defined ex ante as administrative parameters rather than emerging endogenously from system-wide optimization.

\section{Case Study}

\subsection{Brazilian Power System 2034} \label{brazil_case_study}

We simulate the Brazilian Power System for the 2034 planning horizon, representing its four main regions: North (N), Northeast (NE), Southeast (SE), and South (S), with regional generation, demand and transmission constraints.

\textbf{Load characteristics:} Average system demand is approximately 117~GW, with a peak of 151~GW. Demand growth reflects economic expansion, along with emerging loads from electric vehicles, green hydrogen, industrial electrification and data centers. Together with distributed solar, these trends produce an hourly net-load profile characterized by low midday demand and steep evening ramps.

\textbf{Existing generation:} The current system includes 167~hydroelectric plants, 175~thermal plants, 163~wind farms, 154~solar plants, and 4~small hydro facilities. Variable renewable generation totals approximately 78~GW of installed capacity, creating the flexibility challenges that motivate this study.

\textbf{Expansion candidates:} The candidate set consists of 86~GW of potential investments across 13 technologies and 4 regions. Hydro accounts for 7.2~GW (conventional 4.2, small hydro 3), thermal for 29.5~GW (Gas LNG 16, Gas Peaker 8, Gas Pre-Salt 4, coal 1.5), renewables for 43.3~GW (wind 26.7, solar 11.2, offshore wind 3.4, biomass 2), and flexibility resources for 6~GW (batteries 2, pumped hydro 3.2, demand response 0.8).

\subsection{Methodological Framework}

The model employs a two-stage stochastic formulation solved as an equivalent deterministic problem using OptGen~\cite{OptGen}, the complete mathematical formulation is provided in the Appendix. The model uses typical-day modelling, with 1 representative day per month, to capture both daily flexibility requirements and intra-annual seasonality. 
To reduce computational effort without losing the main characteristics of the original sample, the 200 in-sample scenarios were grouped into 20 clusters used in all auction simulations. The resulting in-sample mixed-integer linear program contains approximately 14~million constraints, 15~million variables, and 467~integer investment decisions. One OEA auction instance was solved in 167~minutes on AWS cloud with 50~parallel threads at a 3\% optimality gap tolerance using Xpress~9.1. 
At the expansion horizon considered in this work, the operational model abstracts from unit commitment constraints and detailed nodal power balance. These elements are essential for short-term scheduling problems, but their explicit inclusion in a stochastic, multi-temporal expansion framework would significantly increase problem dimensionality. The adopted representation preserves the cost, flexibility, and reliability trade-offs that drive investment decisions.

\subsection{Four Procurement Approaches}

The experimental design compares four procurement approaches using the same set of project bids. The three credit-based cases are constructed directly from the auction formulation presented in Section \ref{GF_Brazil}. In each case, projects are selected according to Eq. \eqref{Auc_GF}, with the adequacy contribution $\phi_j$ defined either by Physical Guarantee or Capacity Credits and the requirement $\hat{d}$ adjusted accordingly. The sequential case applies these two auctions consecutively.
The suffix “+T” indicates that transmission expansion is optimized after generation procurement. Once the generation portfolio is fixed, the expansion model is rerun with transmission investments as the only remaining decision variables. The four procurement approaches are defined as follows:

\begin{enumerate}
    \item \textbf{Integrated auction (OEA)}: Co-optimizes generation and transmission through a unified expansion model using OptGen~\cite{OptGen}. Investment selection explicitly accounts for operational synergies, locational value, and network constraints.

    \item  \textbf{Energy auction (NEA+T)}: Replicates current energy auctions, selecting projects by their Physical-Guarantee credits, computed for all candidate projects including batteries and demand response.

    \item  \textbf{Capacity auction (CA+T)}: Selects projects by their capacity credits: thermal and hydro plants are procured for 95\% of the auction requirement, and BESS projects for the remaining 5\%. Bids are adjusted by the multiplicative Factor~A for ramp capability, and batteries located in the Northeast receive a 10\% locational discount~($\beta$).
    
    \item  \textbf{Sequential auction (NEA+CA+T)}: First runs the energy auction to meet energy adequacy, then evaluates the residual capacity gap and runs the capacity auction to procure additional firm capacity. 
\end{enumerate}

\begin{figure*}[htbp]
\centering
\includegraphics[width=0.85\textwidth]{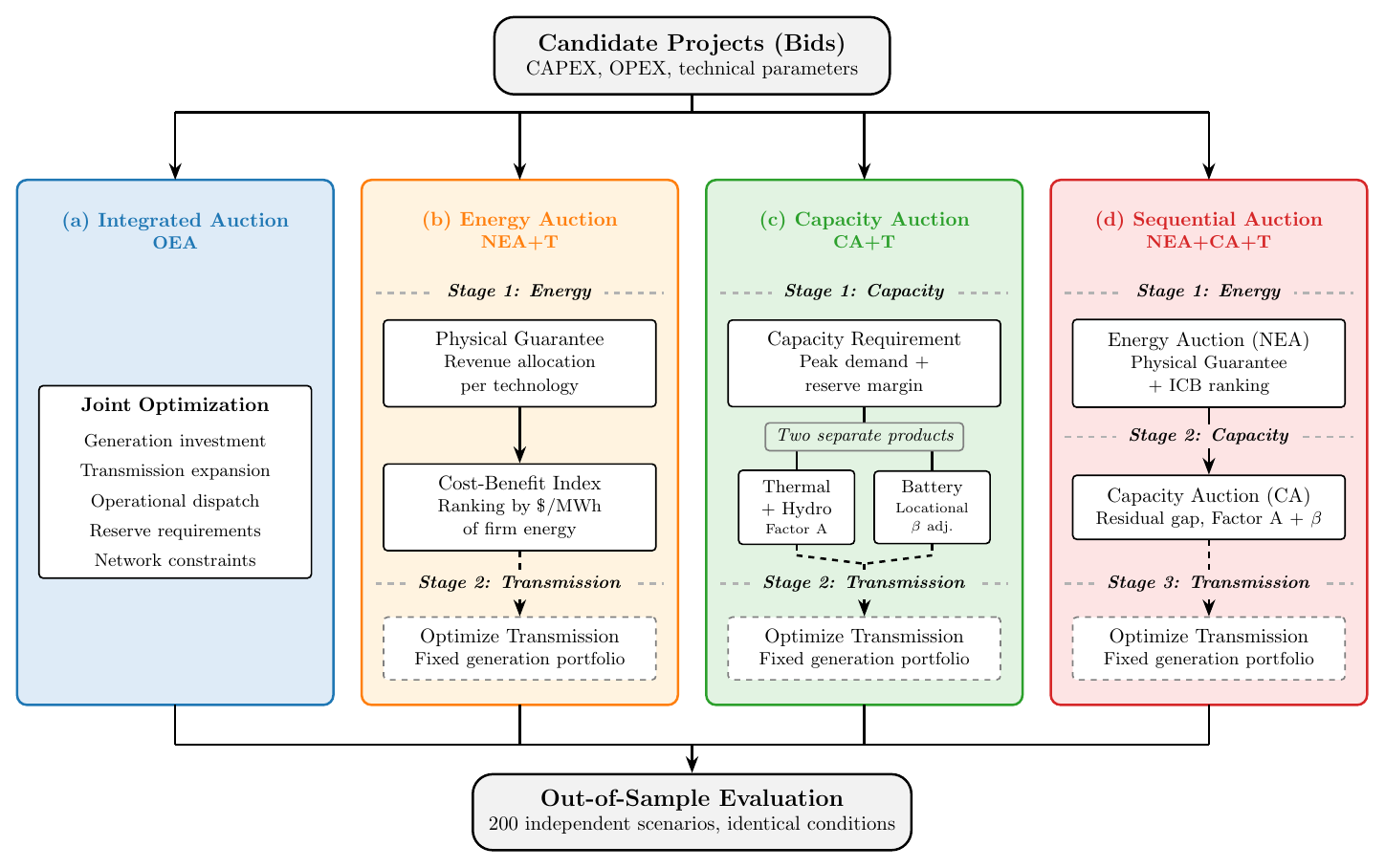}
\caption{Overview of the four procurement approaches. OEA co-optimizes all decisions simultaneously, while the alternative approaches decompose generation and transmission planning into sequential stages (dashed arrows).}
\label{fig:procurement_flowchart}
\end{figure*}

In a real auction, the fixed-revenue bids are strategic and would differ across agents under competition. Recalling the incentive-compatibility concern raised above, we deliberately hold them fixed to isolate the effect of mechanism design from bidding behavior. To compare the four procurement approaches on equal terms, every candidate enters each approach with the same technology-typical annualized investment and fixed O\&M costs as its fixed-revenue requirement, and the same declared variable cost (CVU). Any difference in outcomes is therefore attributable to how each mechanism selects and values resources, not to heterogeneity in bidding behavior. 
Following investment selection under each approach, we conduct identical out-of-sample operational simulations across 200 scenarios drawn independently of those used for selection, to ensure a fair comparison and guard against in-sample overfitting. All portfolios face the same demand realizations and renewable variability, and the cost metrics in Section~IV reflect out-of-sample performance.

\section{Results}

\subsection{System Costs and Emissions}

Table~\ref{tab:cost_comparison} summarizes the cost and emissions outcomes across all four procurement approaches.\footnote{All monetary values are in billions of US dollars, converted from Brazilian Reais at R\$5/US\$, unless otherwise noted. Investment costs are annualized, and operating costs cover the modeled year.} We report both expected cost and CVaR$_{95\%}$ (Conditional Value-at-Risk), which measures the expected operating cost in the worst 5\% of scenarios, a standard risk metric for capturing tail exposure. The integrated approach achieves the lowest total system cost of \$10.93~billion, representing savings of 8\% relative to the energy auction (\$11.82~billion), 16\% relative to the sequential auction (\$12.95~billion), and 21\% relative to the capacity auction (\$13.85~billion).

\begin{table}[htbp]
\centering
\caption{Summary of System Costs and Emissions by Procurement Approach}
\label{tab:cost_comparison}
\setlength{\tabcolsep}{4pt}
\begin{tabular}{@{}lrrrr@{}}
\toprule
\textbf{Metric} & \textbf{OEA} & \textbf{NEA+T} &
\textbf{CA+T} & \textbf{NEA+CA+T} \\
\midrule
Total Expected Cost (\$~bi) & 10.93 & 11.82 & 13.85 & 12.95 \\
Total CVaR$_{95\%}$ Cost (\$~bi) & 24.42 & 29.47 & 30.95 & 27.13 \\
\midrule
Investment Cost (\$~bi) & 4.62 & 2.96 & 4.69 & 4.75 \\
Expected Operating Cost (\$~bi) & 6.32 & 8.86 & 9.16 & 8.20 \\
CVaR$_{95\%}$ Operating Cost (\$~bi) & 19.80 & 26.51 & 26.26 & 22.38 \\
\midrule
Expected CO$_2$ (Mt/yr) & 51.0 & 65.3 & 61.6 & 63.3 \\
\bottomrule
\end{tabular}
\end{table}

The cost decomposition reveals different investment strategies. The energy auction minimizes capital expenditure (\$2.96~billion) but incurs operating costs of \$8.86~billion, second only to the capacity auction (\$9.16~billion), reflecting systematic under-investment that forces continued reliance on existing generation. This under-investment also creates elevated tail risk: the CVaR$_{95\%}$ for the energy auction reaches \$26.51~billion, approximately 34\% higher than the integrated approach, indicating vulnerability to adverse scenarios. Conversely, the capacity and sequential auctions invest similar amounts to the integrated approach (\$4.69~billion and \$4.75~billion, respectively), predominantly in gas-fired technologies, but without achieving proportional reductions in operating costs. The integrated approach balances investment (\$4.62~billion) with operational efficiency (\$6.32~billion), achieving both the lowest expected cost and the lowest tail risk. 

Emissions exhibit a consistent pattern across procurement approaches, although CO$_2$ reduction is not an objective of the optimization. The model minimizes total system cost and does not include a carbon constraint or price. The observed pattern reflects the correlation between fossil emission and thermal variable costs, so cost minimization naturally selects portfolios with lower emissions when cheaper low-carbon alternatives (wind, solar, hydro) are available. The integrated approach produces expected emissions of 51.0~MtCO$_2$ per year, approximately 20\% below the 61.6--65.3~MtCO$_2$ range observed under decoupled procurement. This is a co-benefit of cost-optimal selection, not a direct policy outcome of the mechanism.

\subsection{Selected Portfolios}

Table~\ref{tab:investment_tech} presents the selected portfolios by technology, region, and auction design.

\begin{table}[htbp]
\centering
\caption{Selection by Technology and Region (MW)}
\label{tab:investment_tech}
\footnotesize
\begin{tabular}{@{}llrrrrr@{}}
\toprule
\textbf{Auction} & \textbf{Technology} & \textbf{NE} & \textbf{N} & \textbf{SE} & \textbf{S} & \textbf{Total} \\ \midrule
\multirow{6}{*}{OEA}
 & Wind & 7,400 & 200 & --- & 1,600 & 9,200 \\
 & Hydro & --- & 1,716 & 630 & 982 & 3,327 \\
 & Small Hydro & --- & 1,000 & --- & --- & 1,000 \\
 & Gas Pre-Salt & --- & --- & 4,000 & --- & 4,000 \\
 & Battery & 150 & --- & --- & --- & 150 \\
 & DR & 200 & 200 & 200 & 200 & 800 \\
\midrule
\multirow{3}{*}{NEA+T}
 & Hydro & --- & 774 & 972 & 982 & 2,727 \\
 & Gas Pre-Salt & --- & --- & 4,000 & --- & 4,000 \\
 & Gas LNG & --- & --- & 2,400 & 1,600 & 4,000 \\
\midrule
\multirow{5}{*}{CA+T}
 & Hydro & --- & 1,716 & 972 & 1,363 & 4,050 \\
 & Gas LNG & 1,400 & 1,200 & 1,200 & 1,200 & 5,000 \\
 & Gas Peaker & 2,000 & 2,000 & 2,000 & 2,000 & 8,000 \\
 & Battery & 450 & 200 & 150 & 150 & 950 \\
 & DR & 200 & 200 & 200 & 200 & 800 \\
\midrule
\multirow{6}{*}{\shortstack[l]{NEA+\\CA+T}}
 & Hydro & --- & 1,716 & 972 & 1,363 & 4,050 \\
 & Gas Pre-Salt & --- & --- & 4,000 & --- & 4,000 \\
 & Gas LNG & --- & --- & 2,400 & 1,600 & 4,000 \\
 & Gas Peaker & 1,400 & 1,400 & 1,400 & 1,400 & 5,600 \\
 & Battery & 300 & 50 & 50 & 50 & 450 \\
 & DR & 200 & 200 & 200 & 200 & 800 \\
\bottomrule
\multicolumn{7}{@{}l}{\footnotesize DR: Demand Response.}
\end{tabular}
\end{table}

The integrated portfolio emphasizes diversification and regional complementarity: wind capacity concentrates in the Northeast (7.4~GW), where resource quality is highest, conventional and small hydro expand in the North (2.7~GW), Gas Pre-Salt is located in the Southeast near load centers (4~GW), Demand Response is contracted in all regions, and battery storage provides localized flexibility in the Northeast. The base system already includes the interconnection expansions contracted through previous transmission auctions and scheduled for operation by 2034. Additional reinforcements were offered as candidates, but none were selected, indicating that the already-contracted network expansion was sufficient for the resulting generation portfolio.

Decoupled approaches yield markedly different outcomes. The energy auction selects no variable renewable generation. Instead, it invests heavily in gas (4~GW LNG, 4~GW Pre-Salt) and conventional hydro. The capacity auction favors the selection of gas-fired resources (13~GW of LNG and peaker capacity) distributed uniformly across all four regions, ignoring locational value. The sequential auction combines both patterns with gas portfolios totaling nearly 14~GW. These patterns follow from how the credits are calculated. Wind misses the energy auction narrowly, its cheapest project at \$56/MWh against a clearing marginal cost near \$40/MWh, while storage is excluded by a wide margin, above \$2{,}200/MWh. Under the capacity credit, batteries are contracted through a specific quota, while Gas Pre-Salt is not selected, as its high minimum generation requirement leaves it uncompetitive for the factor A adjustment.

\subsection{Price Signals}

\begin{figure*}[htbp]
    \centering
    \includegraphics[width=0.85\textwidth]{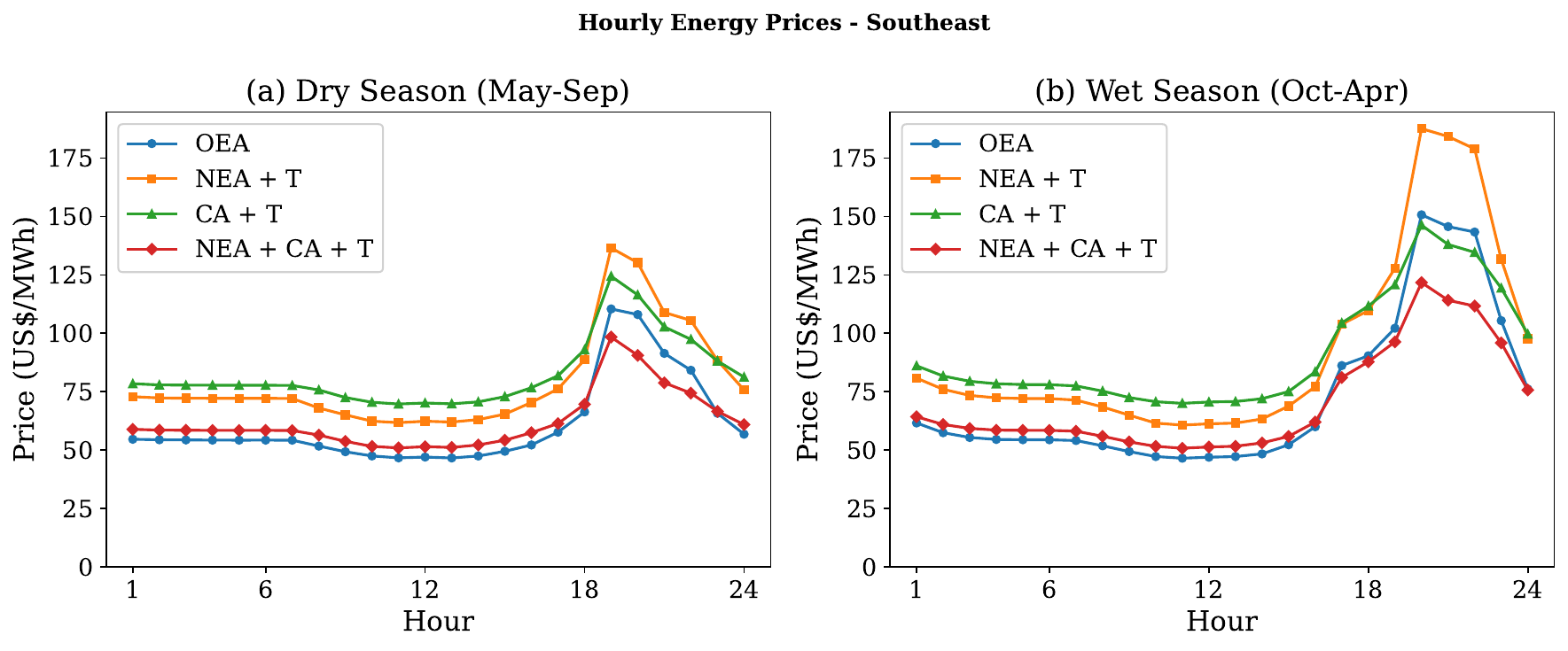}
    \caption{Hourly energy prices in Southeast by season. Shaded areas represent the interquartile range across scenarios.}
    \label{fig:energy_prices_hourly}
\end{figure*}

Given the investments selected under each auction design, Figure~\ref{fig:energy_prices_hourly} presents out-of-sample simulated hourly energy prices (shadow prices on the power balance constraint) in the Southeast subsystem, Brazil’s primary load center, for the dry and wet seasons.

The energy auction exhibits elevated price volatility, with evening peaks reaching \$137/MWh during the dry season and \$188/MWh during the wet season. This pattern reflects the under-invested portfolio's greater reliance on high-cost generation during demand peaks. 

Under the integrated approach, daytime prices are the lowest across all approaches ($\approx$\$47/MWh), reflecting abundant renewable generation in the diversified portfolio. Evening prices then rise sharply as solar output declines, reaching levels comparable to the energy auction: these peaks reflect the true scarcity value of flexibility during the solar-to-evening transition, providing the investment signals needed for marginal resources to recover their fixed costs. By contrast, the capacity and sequential auctions exhibit flatter profiles. The capacity auction sustains the highest overnight prices, \$78/MWh in both seasons, against \$54/MWh under the integrated approach, reflecting its gas-heavy portfolio, while the sequential auction attenuates the evening peak, \$98/MWh dry and \$122/MWh wet, against \$110/MWh and \$151/MWh under the integrated approach. 

\begin{figure*}[htbp]
    \centering
    \includegraphics[width=0.85\textwidth]{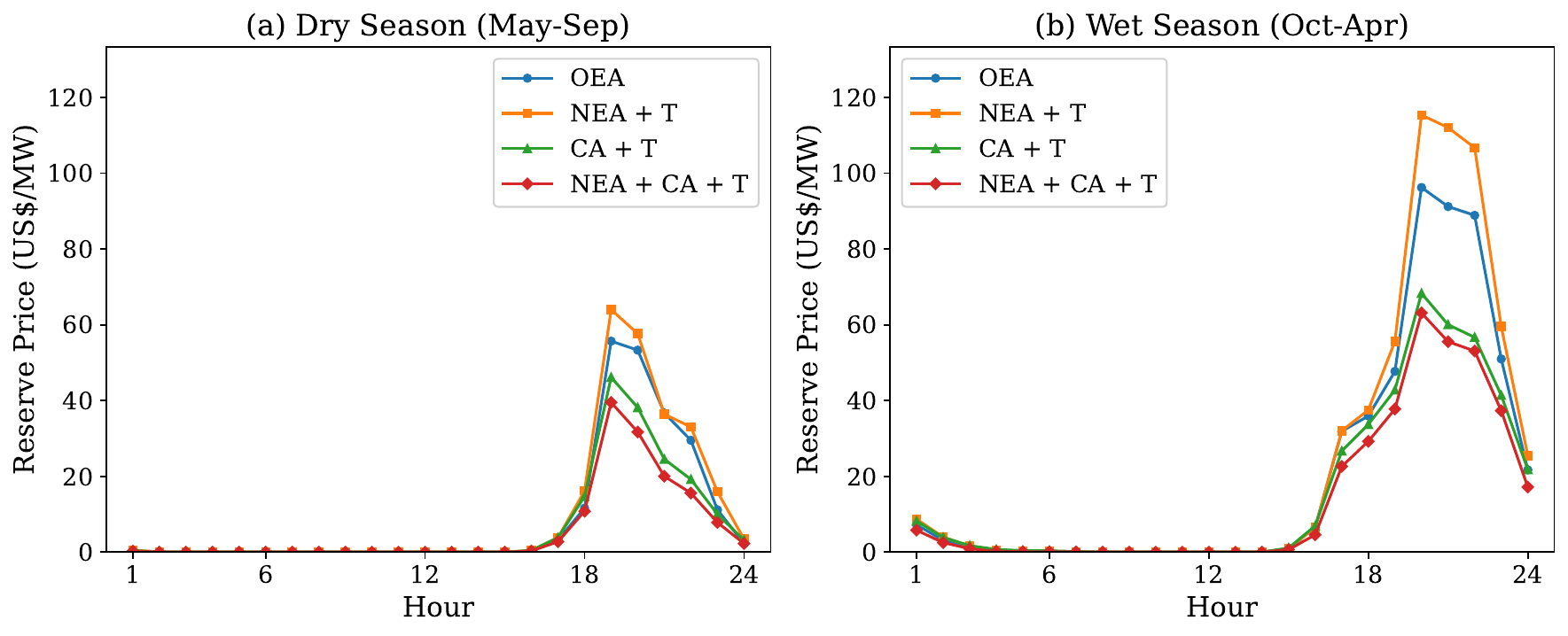}
    \caption{Hourly reserve prices by season. Shaded areas represent the interquartile range across scenarios.}
    \label{fig:reserve_prices_seasonal}
\end{figure*}

\subsection{Reserve Prices}

Figure~\ref{fig:reserve_prices_seasonal} presents hourly reserve prices by season. Reserve prices reflect the marginal cost of providing flexibility to accommodate renewable variability. All approaches exhibit pronounced evening peaks (hours 17--22), when solar generation declines and demand remains elevated. The energy auction shows the highest reserve prices, peaking at \$64/MW in the dry season and \$115/MW in the wet season, reflecting under-investment in flexible capacity. The integrated approach achieves moderate peaks of \$56/MW (dry) and \$96/MW (wet), while the capacity and sequential auctions exhibit the lowest reserve prices (\$40--68/MW). The effect is indirect: their large peaker fleets cover the evening energy requirement, leaving hydro and battery capacity free to supply reserve.

The seasonal pattern reveals higher energy and reserve prices during the wet season across all approaches. Reservoirs enter the wet season at relatively low storage levels after being drawn down during the dry season, increasing water values and scarcity risk until inflows rebuild storage. Reserve prices follow the same pattern because hydro units provide most reserves, creating a trade-off between reserve provision and energy generation.

\subsection{Reliability Analysis}

Table~\ref{tab:reliability} presents reliability metrics across procurement approaches. We report the average deficit in MW (expected unserved power across all hours and scenarios) and the CVaR$_{1\%}$ deficit (expected unserved power in the worst 1\% of scenarios).

\begin{table}[htbp]
\centering
\caption{Reliability Metrics by Approach}
\label{tab:reliability}
\begin{tabular}{@{}lrr@{}}
\toprule
\textbf{Approach} & \textbf{Avg Deficit (MW)} & \textbf{CVaR$_{1\%}$ (MW)} \\ \midrule
OEA          & 2.4   & 88    \\
NEA+T        & 9.8   & 574   \\
CA+T      & 0.02  & 2.0   \\
NEA+CA+T  & 0.02  & 2.2   \\
\bottomrule
\end{tabular}
\end{table}

The results reveal a fundamental trade-off between reliability and economic efficiency. The capacity and sequential auctions achieve near-zero expected deficit (0.02~MW) and minimal tail risk (CVaR$_{1\%}$~$\approx$~2~MW), but at a substantial cost premium driven almost entirely by operating cost (\$9.16~billion and \$8.20~billion versus \$6.32~billion). The energy auction, conversely, under-invests and exhibits the worst reliability result: average deficit reaches 9.8~MW (4$\times$ higher than the integrated approach), but more critically, its tail risk (CVaR$_{1\%}$) reaches 574~MW, over 6 times worse than the integrated approach.

The integrated approach finds an efficient middle ground. Its CVaR$_{1\%}$ deficit of 88~MW represents only 0.07\% of average demand (117~GW), well below the CVaR$_{1\%}$(ENS) $\leq$ 5\% of demand threshold~\cite{LRCAP_EPE_2023}. Rather than over-investing to eliminate the deficit, it permits an economically optimal level of scarcity during peak hours, achieving the lowest total system cost while satisfying the reliability metric.

\section{Conclusions}

This study compared integrated versus decoupled procurement of energy and capacity for power system expansion using a detailed model of the Brazilian system for 2034. The proposed Optimized Expansion Auction achieves 8--21\% lower total system cost relative to decoupled alternatives, with expected costs of \$10.93~billion compared to \$11.82--13.85~billion under the energy, capacity, and sequential auctions. These gains arise from proper valuation of technology synergies that decoupled approaches, which select projects on single-attribute credits, firm energy or firm capacity, systematically overlook.

The proposed approach also achieves the lowest tail cost (CVaR$_{95\%}$ of operating cost, \$19.80~billion versus \$22.38--26.51~billion for decoupled approaches). Since the clearing problem minimizes a risk-neutral expectation, this is a property of the selected portfolio rather than of the mechanism's objective, although the formulation extends naturally to a risk-averse objective, for instance by adding a CVaR term to the clearing problem. As a co-benefit, the proposed approach reduces emissions by approximately 20\%. This reflects the correlation between cost and emission of fossil resources in the current Brazilian candidate set: minimizing total system cost naturally selects portfolios with lower emissions. Explicit carbon policy objectives could also be incorporated into the model through dedicated constraints. In this work, however, they were not included in the formulation and are left for future research.

The investment portfolios reveal fundamentally different strategies. The proposed auction produces diversified portfolios emphasizing wind (9.2~GW), hydro (4.3~GW), and Gas Pre-Salt (4~GW), while decoupled approaches converge on gas-fired technologies---up to 13~GW of LNG and peakers under the capacity auction, that satisfy individual procurement criteria but prove suboptimal when evaluated holistically.  Ultimately, the performance of the OEA relies directly on the fidelity of the underlying mathematical expansion model. As computational power and operational representations advance, the auction's efficiency will automatically advance with it.

Beyond the improvements in project selection, the implementation of the OEA carries profound regulatory implications for any market in which it is adopted. This auction design has the advantage of not requiring administratively determined resource adequacy accreditation metrics, such as the Effective Load Carrying Capability (ELCC) used in U.S. capacity markets or the Physical Guarantee used in Brazil. 
Under the OEA framework, the centralized optimization model endogenously resolves system security and capacity requirements, selecting the units that ensure the reliability needs at least cost. 

Several extensions of the proposed framework could be explored in future research: (1) the design of a joint auction for generation, transmission, and large loads, thereby creating a transparent and efficient framework for selecting the optimal expansion required to meet the growing demand from AI data centers; (2) extending the framework to internalize optimal gas infrastructure expansion within the selected portfolio; (3) analyzing the contracting strategies that may arise from the resulting portfolio of winning resources; and (4) incorporating stability requirements, such as inertia and reactive power control. 

\section*{Acknowledgments}\footnotetext{The case study input data and result files will be made publicly available in an online repository.}
\noindent The authors express their gratitude to Dr. Mario Veiga and Dr. Sérgio Granville for their insightful inputs and inspiration. We also thank PSR for granting access to the PSRCloud platform and for supporting the databases used in Case Study~\ref{brazil_case_study}.

\bibliographystyle{IEEEtran}
\bibliography{main}

\appendices

\section{Detailed Mathematical Formulation} \label{mathematical_formulation}

This appendix presents the complete mathematical formulation underlying the proposed auction mechanism. The model extends the conceptual framework in Section~II with full operational detail.

The expansion model minimizes investment plus expected operational cost. Investment decisions $x_a$ are integer, and operational variables are indexed by scenario $s$, subsystem $n$, and time $t$. The time index is structured hierarchically: months, typical days (1 typical day) within each month, and hourly blocks within each day, capturing both seasonal hydrology and intra-day flexibility. The availability indicator $\chi_a = 1$ for existing assets and $\chi_a = x_a$ for candidates.

{\small
\begin{align}
\min_{x,g,e,v,f,\delta} \; & \sum_{a} c^{inv}_a x_a + \sum_{s} p_s \sum_{t} \Bigl( \sum_{j} c_j g_{s,j,t} + c^{shed} \sum_{n} \delta_{s,n,t} \Bigr) \\
\text{s.t.} \; & \sum_{j \in G_n} g_{s,j,t} + \sum_{h \in H_n} e_{s,h,t} + \sum_{r \in R_n} g^r_{s,r,t} \nonumber \\
& \quad + \sum_{b \in B_n} \left( p^{dis}_{s,b,t} - p^{ch}_{s,b,t} \right) + \sum_{l \in L^{in}_n} f_{s,l,t} \nonumber \\
& \quad - \sum_{l \in L^{out}_n} f_{s,l,t} + \delta_{s,n,t} = Dem_{s,n,t} \\
& e_{s,h,t} = \rho_h \cdot turb_{s,h,t} \\
& v_{s,h,t+1} = v_{s,h,t} + a_{s,h,t} - turb_{s,h,t} - spill_{s,h,t} \nonumber \\
& \quad + \sum_{m \in \Omega_h} (turb_{s,m,t} + spill_{s,m,t}) \\
& v_{s,h,T} = v_{s,h,1} \\
& soc_{s,b,t+1} = soc_{s,b,t} + \eta^{ch}_b p^{ch}_{s,b,t} - p^{dis}_{s,b,t}/\eta^{dis}_b \\
& \sum_{h} r^{hy}_{s,h,t} + \sum_{b} r^{bat}_{s,b,t} \geq \sum_{n} \Bigl( \beta^{d}_{n} Dem_{s,n,t} \nonumber \\
& \qquad + \beta^{w}_{n} \!\!\sum_{r \in \mathcal{W}_n}\!\! g^{r}_{s,r,t} + \beta^{s}_{n} \!\!\sum_{r \in \mathcal{S}_n}\!\! g^{r}_{s,r,t} \Bigr) \label{eq:reserve} \\
& 0 \leq g_{s,j,t} \leq \overline{G}_j \chi_j, \quad |g_{s,j,t} - g_{s,j,t-1}| \leq \Delta_j \\
& 0 \leq e_{s,h,t} \leq \overline{E}_h \chi_h, \quad \underline{V}_h \chi_h \leq v_{s,h,t} \leq \overline{V}_h \chi_h \\
& 0 \leq g^r_{s,r,t} \leq \hat{g}_{s,r,t} \chi_r, \quad |f_{s,l,t}| \leq \overline{F}_l \chi_l \\
& \underline{S}_{b} \chi_{b} \leq soc_{s,b,t} \leq \overline{S}_{b} \chi_{b}, \; 0 \leq p^{ch}_{s,b,t}, p^{dis}_{s,b,t} \leq \overline{P}_{b} \chi_{b} \\
& 0 \leq r^{hy}_{s,h,t} \leq \overline{E}_h \chi_h - e_{s,h,t} \nonumber \\
& 0 \leq r^{bat}_{s,b,t} \leq \overline{P}_{b} \chi_{b} - p^{dis}_{s,b,t}
\end{align}
}

\noindent\textbf{Sets:}
{\footnotesize
$G_n$, $H_n$, $R_n$, $B_n$: thermal, hydro, renewable, and battery units in subsystem $n$;
$\mathcal{W}_n$, $\mathcal{S}_n$: wind and solar plants in subsystem $n$ (subsets of $R_n$);
$L^{in}_n$, $L^{out}_n$: interconnections into and out of subsystem $n$;
$\Omega_h$: hydro plants immediately upstream of $h$.
}

\noindent\textbf{Parameters:}
{\footnotesize
$c^{inv}_a$: annualized investment cost (the submitted fixed-revenue bid, set here to each technology's typical value);
$c_j$, $c^{shed}$: thermal variable cost, value of lost load;
$p_s$: scenario probability;
$\overline{G}_j$, $\Delta_j$: thermal capacity, ramp limit;
$\rho_h$, $a_{s,h,t}$: hydro productivity, inflow;
$\underline{V}_h$, $\overline{V}_h$, $\overline{E}_h$: reservoir bounds, max generation;
$\hat{g}_{s,r,t}$: renewable availability;
$\eta^{ch}_b$, $\eta^{dis}_b$, $\underline{S}_{b}$, $\overline{S}_{b}$, $\overline{P}_{b}$: battery efficiencies, energy and power bounds;
$\overline{F}_l$: interconnection capacity;
$\beta^{d}_n$, $\beta^{w}_n$, $\beta^{s}_n$: per-subsystem reserve coefficients.
}

\noindent\textbf{Variables:}
{\footnotesize
$x_a$: investment decision (units);
$g_{s,j,t}$, $e_{s,h,t}$, $g^r_{s,r,t}$: thermal, hydro, renewable generation;
$v_{s,h,t}$, $turb_{s,h,t}$, $spill_{s,h,t}$: reservoir volume, turbined and spilled water;
$p^{ch}_{s,b,t}$, $p^{dis}_{s,b,t}$, $soc_{s,b,t}$: battery charge, discharge, SOC;
$f_{s,l,t}$, $\delta_{s,n,t}$: interconnection flow, load shedding;
$r^{hy}_{s,h,t}$, $r^{bat}_{s,b,t}$: reserve provision.
}

\end{document}